\documentclass[12pt,onecolumn]{IEEEtran}

\renewcommand{\baselinestretch}{1.3}
\usepackage{cite}
\usepackage{subcaption}
\usepackage[font={small}]{caption}
\usepackage{algorithm}
\usepackage{algorithmic}
\usepackage[algo2e]{algorithm2e} 

\usepackage{mathtools}
\DeclarePairedDelimiter\ceil{\lceil}{\rceil}
\DeclarePairedDelimiter\floor{\lfloor}{\rfloor}

\usepackage{array}
\newcolumntype{P}[1]{>{\centering\arraybackslash}p{#1}}

\usepackage{color}
\usepackage{subcaption}
\usepackage{graphicx}
\usepackage{amsmath}
\usepackage{hyperref}
\usepackage{amssymb}
\usepackage{amsthm}
\usepackage{dsfont}

\newtheorem*{proof*}{\hspace{0.8cm}{\it{Proof}}}

\newcounter{relctr} 
\everydisplay\expandafter{\the\everydisplay\setcounter{relctr}{0}} 

\newcommand\labelrel[2]{%
  \begingroup
    \refstepcounter{relctr}%
    \stackrel{\textnormal{(\alph{relctr})}}{\mathstrut{#1}}%
    \originallabel{#2}%
  \endgroup
}
\AtBeginDocument{\let\originallabel\label} 

\begin{document}

\title{\LARGE  No Cross-Validation Required: An Analytical Framework for Regularized Mixed-Integer Problems (Extended Version)$^{\dagger}$ \thanks{$\dagger$ {This work is the extended version of \cite{Letter}, including all proofs and more simulation studies.}} }

\author{\small
\hspace{-6pt}\begin{tabular}[t]{c@{\extracolsep{5em}}c@{\extracolsep{5em}}c} Behrad~Soleimani  & Behzad~Khamidehi & Maryam~Sabbaghian\\
\textit{Department of ECE} & \textit{Department of ECE}& \textit{Department of ECE} \\
\textit{University of Maryland}  & \textit{University of Toronto} & \textit{University of Tehran}\\
College Park, MD, USA & Toronto, ON, Canada & Tehran, Iran \\
behrad@umd.edu & b.khamidehi@utoronto.ca & msabbaghian@ut.ac.ir\\
\end{tabular}
\vspace{-10mm}
}
	\maketitle
	\begin{abstract}
	This paper develops a method to obtain the optimal value for the regularization coefficient in a general mixed-integer problem (MIP). This approach eliminates the cross-validation performed in the existing penalty techniques to obtain a proper value for the regularization coefficient. We obtain this goal by proposing an alternating method to solve MIPs. First, via regularization, we convert the MIP into a more mathematically tractable form. Then, we develop an iterative algorithm to update the solution along with the regularization (penalty) coefficient. We show that our update procedure guarantees the convergence of the algorithm. Moreover, assuming the objective function is continuously differentiable, we derive the convergence rate, a lower bound on the value of regularization coefficient, and an upper bound on the number of iterations required for the convergence. We use a radio access technology (RAT) selection problem in a heterogeneous network to benchmark the performance of our method. Simulation results demonstrate near-optimality of the solution and consistency of the convergence behavior with obtained theoretical bounds.
	\end{abstract}
	\vspace{-5mm}
	\IEEEpeerreviewmaketitle
\section{Introduction}
Mixed-integer problems (MIPs) include a variety of problems, most of which are broadly applied in telecommunications systems. This includes assignment or classification of variables whose specific application in telecom can be resource allocation \cite{boyd2004convex}. A conventional technique to solve MIPs is to convert the problem into a continuous form, relax the integer variables, and regularize (penalize) the objective function. The regularization function forces the corresponding relaxed variables to be integers. Otherwise, the objective function is significantly penalized. To control this penalty, the regularization function is scaled by a regularization coefficient. This coefficient highly impacts the results. Hence, finding its optimal value is of vital importance \cite{ModelSelection}.

In most cases, there is not an analytical closed-form expression for the optimal regularization coefficient. Thus, the optimal coefficient has to be obtained by the time consuming process of cross-validation. For instance, in \cite{DCref}, \cite{ng2015secure}, \cite{KTH}, \cite{BehzadDC}, and \cite{BehradCluster} where a quadratic regularization function is used, the value of the regularization coefficients is attained via cross-validation. This limits the applications of the approach as it cannot be applied to high-dimensional problems due to its increasing computational complexity.

\begin{figure}
        \centering
        \includegraphics[trim={10cm 7cm 10cm 1cm},width=0.5\linewidth]{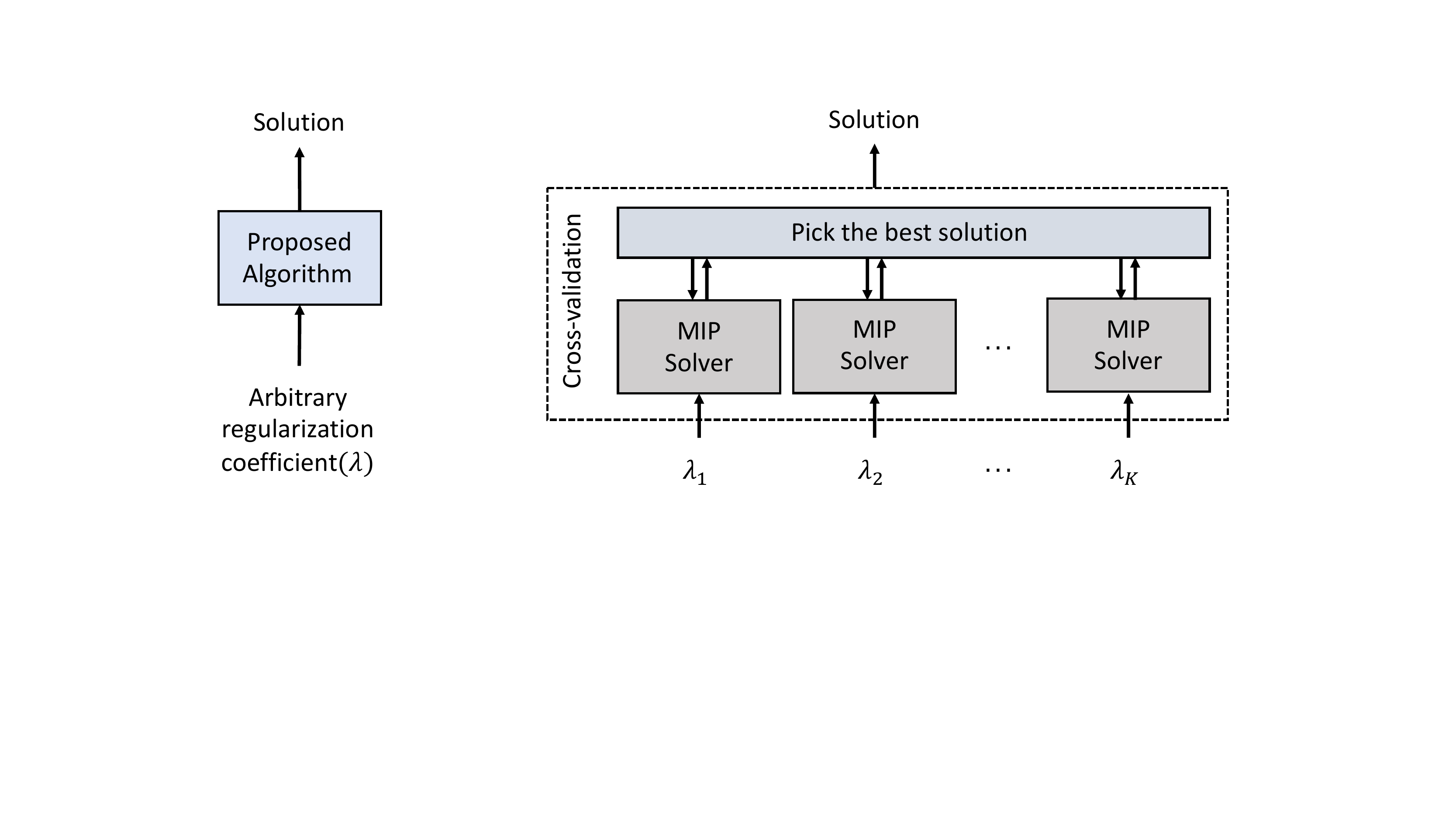}
        \caption{Overview of cross-validation procedure and the proposed algorithm.}
        \label{fig:Intro}
    \end{figure}

In this paper, in a mathematical framework, we obtain a closed-form expression for the optimal value of the regularization coefficient. Thus, we avoid the extensive complexity involved in cross-validation as depicted in Fig. \ref{fig:Intro}. Through theorems that we prove, we also derive the convergence rate and an upper bound on the number of required iterations.  Our proposed alternating algorithm to solve MIPs is as follows. First, we regularize the objective function to reformulate the problem as an equivalent continuous form. Then, to eliminate cross-validation, we utilize an iterative procedure whose solution and regularization coefficient are updated simultaneously. This approach guarantees to achieve a sequence of improved solutions. In fact, we consider MIPs with continuously differentiable objective functions and obtain a lower bound for the regularization coefficient.

To investigate the performance of our proposed method, we apply it to a radio access technology (RAT) selection problem in heterogeneous networks (HetNets) \cite{keshavarz2016hetnets,RATBehrad}. Simulation results show that our method achieves near-optimal solution and the convergence behavior matches the theoretical bounds. 

\section{Preliminaries and Problem Definition}

In this section, we describe the notation used throughout the paper along with the considered optimization problem.\\
    \noindent {\it{Notation}}: We use bold letters to denote vectors. The $i$-th element of vector $\mathbf{x}$ and its $\ell_2$ norm are presented by $x_i$ and $\| \mathbf{x} \|$, respectively. The inner product of two vectors $\mathbf{x}$ and $\mathbf{y}$ is defined as $\langle \mathbf{x} , \mathbf{y} \rangle :=  \sum\limits_{i=1}^{n} x_i y_i$. Function $f(\mathbf{x},\mathbf{y})$ has \textit{Lipschitz} continuity on $\mathbf{x}$ with constant $L$ if
    \begin{align*}
     \left|f(\mathbf{x}_1,\mathbf{y})-f(\mathbf{x}_2,\mathbf{y}) \right| \leq L \left\|\mathbf{x}_1-\mathbf{x}_2 \right\|,
    \end{align*} 
    for all $\mathbf{x}_1$, $\mathbf{x}_2$, and $\mathbf{y}$ on its domain, or equivalently $\left\|\nabla_{\mathbf{x}}f(\mathbf{x},\mathbf{y})\right\| \leq L$ where $\nabla_{\mathbf{x}}$ denotes the gradient (subgradient) operator with respect to $\mathbf{x}$.  {The continuous interval from $0$ to $1$ is presented by $[0,1]$ and the binary set including $0$ and $1$ is shown by $\{0,1\}$.} 
    We also denote the hard-thresholding version of $\mathbf{x} \in [0,1]^n$ with $\widetilde{\mathbf{x}}\in \{0,1\}^n$. The $i$-th element of $\widetilde{\mathbf{x}}$ is
    \begin{align*}
     \widetilde{x}_i = \floor{x_i + \frac{1}{2}},  \ i=1,2,\dots,n,   
    \end{align*}
    where $\floor{.}$ is the floor operator. The distance between $\mathbf{x} \in [0,1]^n$ and $\widetilde{\mathbf{x}}$ is defined as  
    \begin{align*}
    d(\mathbf{x}) := \| \mathbf{x}-\widetilde{\mathbf{x}} \|_1 = n - \sum\limits_{i=1}^{n} \max \{{x}_i , 1-{x}_i \},
    \end{align*}
    where $\| .\|_1$ shows the $\ell_1$ norm. Moreover, function $f(\mathbf{x},\mathbf{y})$ is biconvex in $(\mathbf{x},\mathbf{y})$, if it is convex with respect to each variable while the other one is fixed \cite{gorski2007biconvex}. Similarly, set $B \subseteq X \times Y$ is biconvex on $X \times Y$ if for every fixed $y \in Y$, $B_y = \{x \in X: \ (x,y) \in B \}$ is a convex set in $X$ and for every fixed $x \in X$, $B_x = \{y \in Y: \ (x,y) \in B \}$ is a convex set in $Y$.

\noindent {\it{Problem definition}}: We consider a general MIP with an $m \times 1$ continuous vector variable $\mathbf{y} \in \mathbb{R}^m$ and an $n \times 1$ binary vector variable $\mathbf{x} \in \{ 0,1\}^n$. The optimization problem is given as 
	\begin{align}
	\label{eq:raw-problem}
	\underset{\mathbf{x},\mathbf{y}}{\min} 
	 \ f(\mathbf{x},\mathbf{y})  \ \ \ \
	\text{s.t.}  \ \ \
	\mathbf{x},\mathbf{y} \in \mathcal{D}, \ \ \mathbf{x} \in \{0 , 1\}^n,
	\end{align}
	where $f(\mathbf{x},\mathbf{y}): \ \mathbb{R}^m \times \mathbb{R}^n \rightarrow \mathbb{R}$ is the objective function and $\mathcal{D}$ presents the feasible region 
	\begin{align}
	\label{feasiblereg}
	    \mathcal{D} := \Big\{(\mathbf{x},\mathbf{y}) \Big| &g_i(\mathbf{x},\mathbf{y}) \leq 0, \  i=1,2,\dots,p, \nonumber \\
	    &h_i(\mathbf{x},\mathbf{y}) = 0, \  i=1,2,\dots,q \Big \},
	\end{align}
	such that $g_i(\mathbf{x},\mathbf{y}): \ \mathbb{R}^m \times \mathbb{R}^n \rightarrow \mathbb{R}, \ \forall i$ indicate inequality constraints and $h_i(\mathbf{x},\mathbf{y}): \ \mathbb{R}^m \times \mathbb{R}^n \rightarrow \mathbb{R}, \ \forall i$ are affine functions with respect to both $\mathbf{x}$ and $\mathbf{y}$. {Moreover, we assume that the solution of problem \eqref{eq:raw-problem} over $\mathcal{D}$ is bounded.}
\section{Regularized Equivalent Problem}
    In this section, we reformulate problem \eqref{eq:raw-problem} using a regularization function. The following lemma explains the idea behind this regularization.
    
    \noindent \textbf{Lemma 1.} \textit{ Let $\mathcal{S}$ denote the set 
    \begin{align*}
        \mathcal{S} := \left\{(\mathbf{x},\mathbf{a}) \Big| \langle\mathbf{x},\mathbf{a}\rangle + \langle\mathbf{1}-\mathbf{x},\mathbf{1}-\mathbf{a}\rangle = n , \ \mathbf{x},\mathbf{a} \in [0,1]^n\right\}.
    \end{align*}
    For any pair $(\mathbf{x},\mathbf{a}) \in \mathcal{S}$, we have $\mathbf{x},\mathbf{a} \in \{0,1\}^n$ and $\mathbf{a} = \mathbf{x}$.}
    
    \noindent \textit{Proof.} Let us define vectors $\mathbf{u} = [\mathbf{1}-\mathbf{x} ; \mathbf{x}]$ and $\mathbf{v}=[\mathbf{1}-\mathbf{a};\mathbf{a}]$. Using the Cauchy-Schwarz inequality for these two vectors, we have 
    \begin{multline}
    \label{eq:cauchy}
        \left(\|\mathbf{x} \|^2 +\|\mathbf{1} - \mathbf{x} \|^2  \right)\left( \|\mathbf{a} \|^2 + \|\mathbf{1}-\mathbf{a}\|^2 \right) \geq 
        \left( \langle \mathbf{x},\mathbf{a} \rangle + \langle \mathbf{1}-\mathbf{x}, \mathbf{1}-\mathbf{a} \rangle \right)^2 = n^2, \hspace{0.5cm}\forall (\mathbf{x},\mathbf{a}) \in \mathcal{S}.
    \end{multline}
    On the other hand, since $\mathbf{x} \in [0,1]^n$, we have $\|\mathbf{x} \|^2 +\|\mathbf{1}-\mathbf{x}\|^2 \leq n$, where the equality occurs if and only if $\mathbf{x} \in \{0,1\}^n$. Following the same argument for vector $\mathbf{a}$, we have 
    \begin{align}
    \label{eq:upper}
        \left(\|\mathbf{x} \|^2 +\|\mathbf{1} - \mathbf{x}\|^2  \right)\left( \|\mathbf{a} \|^2 + \|\mathbf{1}-\mathbf{a}\|^2 \right) \leq n^2.
    \end{align}
    Considering \eqref{eq:cauchy} and \eqref{eq:upper}, we can conclude that the equality holds. Thus, $\mathbf{x},\mathbf{a} \in \{0,1\}^n$. Also, using the equality condition of Cauchy-Schwarz inequality, we have $ \mathbf{a}= \mathbf{x} $. \hfill\(\blacksquare\) \\

    Using this lemma and considering auxiliary variable $\mathbf{a}$, the non-convex constraint $\mathbf{x} \in \{0,1 \}^n$ is equivalent to $(\mathbf{x},\mathbf{a}) \in \mathcal{S}$. Set $\mathcal{S}$ is biconvex with respect to  $\mathbf{x}$ and $\mathbf{a}$. Later, in section IV, we explain how biconvexity of $\mathcal{S}$ is used to solve the problem. Now, we can define a regularization function to convert MIP $\eqref{eq:raw-problem}$ into a continuous problem. The following theorem introduces this regularization.
    
    \noindent \textbf{Theorem 1 (Regularized problem).} \textit{The MIP in \eqref{eq:raw-problem} is identical to the following continuous problem for some $\lambda > 0$,}
    \begin{align}
	\label{eq:deductive-penalty}
	&\underset{\mathbf{x},\mathbf{y},\mathbf{a}}{\min} 
	\ \ \ \ \ \ \mathcal{L}_{\lambda} (\mathbf{x},\mathbf{y},\mathbf{a}):=f(\mathbf{x},\mathbf{y}) + \phi_{\lambda}(\mathbf{x},\mathbf{a})  \\
	&\text{s.t.}  
	\ \ \ \ \ \ \ \mathbf{x},\mathbf{y} \in \mathcal{D}, \ \ \mathbf{x},\mathbf{a} \in [0 , 1]^n. \nonumber
	\end{align}
	\textit{where $\phi_{\lambda}(\mathbf{x},\mathbf{a})=\lambda \left(n - \langle\mathbf{x},\mathbf{a}\rangle - \langle\mathbf{1}-\mathbf{x},\mathbf{1}-\mathbf{a}\rangle\right)$.}
	\
	\\
	\\
    \noindent \textit{Proof.}
    {Let $\mathcal{R} := \left\{(\mathbf{x},\mathbf{y},\mathbf{a}) \big |  \mathbf{x},\mathbf{y} \in \mathcal{D}, \mathbf{x},\mathbf{a} \in [0 , 1]^n\right\}$.} Problem \eqref{eq:raw-problem} and its dual are expressed as

    \noindent $\underset{{(\mathbf{x},\mathbf{y},\mathbf{a})\in \mathcal{R}}}{\min} \ \underset{\lambda \geq 0}{\max} \ \mathcal{L}_{\lambda} (\mathbf{x},\mathbf{y},\mathbf{a})$    
    and $\underset{\lambda \geq 0}{\max} \  \underset{{(\mathbf{x},\mathbf{y},\mathbf{a})\in \mathcal{R}}}{\min}  \ \mathcal{L}_{\lambda} (\mathbf{x},\mathbf{y},\mathbf{a})$, respectively. According to the weak duality theorem \cite{boyd2004convex}, we have 
    \begin{equation}
    \label{weak_duality}
    \underset{\lambda \geq 0}{\max} \ \underset{{(\mathbf{x},\mathbf{y},\mathbf{a})\in \mathcal{R}}}{\min}  \  \mathcal{L}_{\lambda} (\mathbf{x},\mathbf{y},\mathbf{a})
    \leq 
    \underset{{(\mathbf{x},\mathbf{y},\mathbf{a})\in \mathcal{R}}}{\min} \  \underset{\lambda \geq 0}{\max} \ \mathcal{L}_{\lambda} (\mathbf{x},\mathbf{y},\mathbf{a}).
    \end{equation}
Defining $\omega(\lambda) := \underset{{(\mathbf{x},\mathbf{y},\mathbf{a})\in \mathcal{R}}}{\min} \ \mathcal{L}_{\lambda} (\mathbf{x},\mathbf{y},\mathbf{a})$, we have 
\begin{align}
\label{Omega_opt}    
\omega(\lambda ^*) = \max_{\lambda\geq 0} \ \omega(\lambda)  
\geq \underset{{(\mathbf{x},\mathbf{y},\mathbf{a})\in \mathcal{R}}}{\min}  \ \mathcal{L}_{\lambda} (\mathbf{x},\mathbf{y},\mathbf{a}).
\end{align}
Based on the Cauchy-schwarz inequality, $\phi_{\lambda}(\mathbf{x},\mathbf{a}) \geq 0$. First, let us consider the case where $\phi_{\lambda}(\mathbf{x},\mathbf{a})=0$. In this case, regardless of the value of $\lambda$, we have $\mathcal{L}_{\lambda} (\mathbf{x},\mathbf{y},\mathbf{a}) = f(\mathbf{x},\mathbf{y})$. In other words, $\mathcal{L}_{\lambda} (\mathbf{x},\mathbf{y},\mathbf{a})$ is constant with respect to $\lambda$. As a result, $\mathcal{L}_{\lambda} (\mathbf{x},\mathbf{y},\mathbf{a}) = \max_{\lambda\geq 0} \ \mathcal{L}_{\lambda} (\mathbf{x},\mathbf{y},\mathbf{a}) $. By substituting this into the right side of \eqref{Omega_opt}, we obtain
\begin{equation}
\label{Proof:case1:1}
\underset{\lambda \geq 0}{\max} \ \underset{{(\mathbf{x},\mathbf{y},\mathbf{a})\in \mathcal{R}}}{\min} \ \mathcal{L}_{\lambda} (\mathbf{x},\mathbf{y},\mathbf{a}) \geq \underset{{(\mathbf{x},\mathbf{y},\mathbf{a})\in \mathcal{R}}}{\min} \  \underset{\lambda \geq 0}{\max}  \ \mathcal{L}_{\lambda} (\mathbf{x},\mathbf{y},\mathbf{a}).
\end{equation}
Comparing \eqref{weak_duality} and \eqref{Proof:case1:1}, we conclude {that the duality gap is zero and strong duality holds as } 
\begin{align}
\label{Proof:case1:result}
\underset{\lambda \geq 0}{\max} \ \underset{{(\mathbf{x},\mathbf{y},\mathbf{a})\in \mathcal{R}}}{\min} \ \mathcal{L}_{\lambda} (\mathbf{x},\mathbf{y},\mathbf{a}) &= \underset{{(\mathbf{x},\mathbf{y},\mathbf{a})\in \mathcal{R}}}{\min} \ \underset{\lambda \geq 0}{\max}  \ \mathcal{L}_{\lambda} (\mathbf{x},\mathbf{y},\mathbf{a})  = \ \omega(\lambda^{\star}).
\end{align}
Due to the fact that function $\omega (\lambda)$ is an increasing function with respect to $\lambda$ \textcolor{black}{(see Appendix \ref{app:omega})}, for $ \lambda > \lambda^{\star}$, we have $\omega(\lambda) \geq \omega(\lambda^{\star})$. However, based on \eqref{Omega_opt}, $\omega(\lambda^{\star})   
\geq \omega(\lambda)$. As a result, $\omega(\lambda^{\star}) = \omega(\lambda) $, $\forall \lambda > \lambda^{\star}$. 
{Combining this with \eqref{Proof:case1:result} yields 
\begin{equation}
\omega(\lambda) = \omega(\lambda^{\star}) = 
\underset{{(\mathbf{x},\mathbf{y},\mathbf{a})\in \mathcal{R}}}{\min} \ \underset{\lambda \geq 0}{\max}\ \mathcal{L}_{\lambda} (\mathbf{x},\mathbf{y},\mathbf{a}), \forall \lambda \geq \lambda^{\star}.
\end{equation}
Since $\underset{{(\mathbf{x},\mathbf{y},\mathbf{a})\in \mathcal{R}}}{\min} \ \underset{\lambda \geq 0}{\max} \ \mathcal{L}_{\lambda} (\mathbf{x},\mathbf{y},\mathbf{a})$ is the solution of \eqref{eq:raw-problem} and $\omega(\lambda)$ is the solution of problem \eqref{eq:deductive-penalty}, we conclude that these two problems are identical and reach the same solutions.
}
Now, we consider the second case where $\phi_{\lambda}(\mathbf{x},\mathbf{a}) > 0$. In this case, since $\omega(\lambda)$ is monotonically increasing with respect to $\lambda$, we have $\omega(\lambda^{\star}) \rightarrow \infty$. This contradicts \eqref{weak_duality} which states that $\omega(\lambda^{\star})$ is finite and upper-bounded by $f^{\star}$, where $f^{\star}$ is the optimal value of \eqref{eq:raw-problem}. Therefore, at the optimal point, $\phi_{\lambda}(\mathbf{x},\mathbf{a}) = 0$, and the results of the first case hold. \hfill\(\blacksquare\) \\

    After this reformulation, in the next section, we derive a method to solve the regularized problem \eqref{eq:deductive-penalty}.
    
    \section{Alternating Algorithm}
    In this section, we develop an alternating algorithm to solve problem \eqref{eq:deductive-penalty}. This method successively minimizes the problem over one variable while fixing all others \cite{boyd2004convex}. {To achieve our goal, first, we consider a biconvex problem and propose our algorithm for this case. We derive some important properties regarding convergence of the algorithm and the update of the regularization parameter. Then, using the obtained properties, we extend our study to a general case and propose an algorithm to solve MIPs which are not necessarily biconvex. }
    \subsection{Biconvex case} 
    {In this section, we derive an algorithm to solve the MIP considering the following assumptions.  
    
    \vspace{0.1cm}
    \noindent {\it{Assumptions:}} The objective function $f(\mathbf{x},\mathbf{y})$ is biconvex in $\mathbf{x}$ and $\mathbf{y}$ and the feasible region $\mathcal{D}$ is convex.} 
    \vspace{0.1cm}

    \noindent Using Theorem 1, we can alter the original MIP into the continuous problem \eqref{eq:deductive-penalty}. Due to biconvexity of $f(\mathbf{x},\mathbf{y})$, $\mathcal{L}_{\lambda} (\mathbf{x},\mathbf{y},\mathbf{a})$ is convex in each argument while others are fixed. Hence, we can use an alternating method to solve problem \eqref{eq:deductive-penalty}. \textcolor{black}{This method is an iterative procedure such that at each iteration, we sequentially solve the problem for one variable among $\mathbf{x},\mathbf{y},\mathbf{a}$ while the others are fixed.} To find a proper value for $\lambda$, we also update its value at each iteration. At the $t$-th iteration, $\lambda^{(t)}=\rho^t\lambda^{(0)}$ where $\rho > 1$ is a known constant and $\lambda^{(0)}$ is the initial value. This update ensures that at some iteration $t$, $\lambda^{(t)} > \lambda^\star$. {This satisfies the required condition in Theorem 1. \textcolor{black}{As a result, there is no need to run cross-validation to obtain $\lambda$.} \textcolor{black}{It is noteworthy to mention that this update works for any arbitrary choice of ($\rho > 1$,$\lambda>0$) (the details are provided in Appendix \ref{app:arbit}).} The overall procedure is presented in Algorithm \ref{alg:biconvex}. It is worth mentioning that the corresponding sub-problem for $\mathbf{a}$ has a closed-form solution. 
   At the $t$-th iteration, the sub-problem is reduced to
    \begin{align}
    \label{eq:update-a}
        \mathbf{a}^{(t)} = \underset{\mathbf{a}}{\text{argmax}} \ \langle\mathbf{a},2\mathbf{x}^{(t)}-\mathbf{1}\rangle \ \ \text{s.t.} \ \ \mathbf{a} \in [0,1]^n,
    \end{align}
    resulting in $\mathbf{a}^{(t)} = \widetilde{\mathbf{x}}^{(t)}$.
    
    \noindent {In the following proposition, we show that Algorithm \ref{alg:biconvex} generates a non-increasing sequence. \\
    \noindent \textbf{\textcolor{black}{Proposition} (Monotone convergence).} 
    \textit{The solution of Algorithm \ref{alg:biconvex} results in a non-increasing sequence. Moreover, this sequence converges to a solution that can not be improved.}\\
    \noindent \textit{Proof.} First, we show that the sequence $\left\{ \mathcal{L}_{\lambda} (\mathbf{x}^{(t)},\mathbf{y}^{(t)},\mathbf{a}^{(t)}) \right\}$ is non-increasing. For any given $t$, we have
    \begin{align*}
        \mathcal{L}_{\lambda} (\mathbf{x}^{(t)},\mathbf{y}^{(t)},\mathbf{a}^{(t)}) \geq         \mathcal{L}_{\lambda} (\mathbf{x}^{(t)},\mathbf{y}^{(t+1)},\mathbf{a}^{(t)}) \geq 
         \mathcal{L}_{\lambda} (\mathbf{x}^{(t+1)},\mathbf{y}^{(t+1)},\mathbf{a}^{(t)}) \geq \mathcal{L}_{\lambda} (\mathbf{x}^{(t+1)},\mathbf{y}^{(t+1)},\mathbf{a}^{(t+1)}),
    \end{align*}
    where the inequalities hold directly because of the updates on variables $\mathbf{x}$, $\mathbf{y}$, and $\mathbf{a}$ in Algorithm \ref{alg:biconvex}. Since problem \eqref{eq:raw-problem} has a finite solution, according to monotone convergence theorem \cite{royden1988real}, the sequence converges to $\underset{t}{\text{inf}} \ \mathcal{L}_{\lambda} (\mathbf{x}^{(t)},\mathbf{y}^{(t)},\mathbf{a}^{(t)})$.}  \hfill\(\blacksquare\) \\
    
       \begin{algorithm}[H]
        \SetAlgoLined
        Set $t = 0$ and $\mathbf{x}^{(t)} = \mathbf{0}$, $\mathbf{a}^{(t)} = \mathbf{0}$.\;\\
        Choose arbitrary $\lambda > 0$, $\rho > 1$.\;\\
        \Repeat{convergence of $\mathbf{x}$, $\mathbf{y}$, and $\mathbf{a}$}{
        Solve the following convex sub-problems {via \textit{interior point}}:\; \\
        $\mathbf{y}^{(t+1)} = \underset{\mathbf{y}}{\text{argmin}} \ \mathcal{L}_{\lambda} (\mathbf{x}^{(t)},\mathbf{y},\mathbf{a}^{(t)}) \ \ \text{s.t.} \ \ \mathbf{y} \in \mathcal{D}$.\; \\
        $\mathbf{x}^{(t+1)} = \underset{\mathbf{x}}{\text{argmin}} \ \mathcal{L}_{\lambda} (\mathbf{x},\mathbf{y}^{(t+1)},\mathbf{a}^{(t)}) \ \ \text{s.t.} \ \ \mathbf{x} \in [0,1]^n \cap \mathcal{D}$. \; \\
        $\mathbf{a}^{(t+1)} = \underset{\mathbf{a}}{\text{argmin}} \ \mathcal{L}_{\lambda} (\mathbf{x}^{(t+1)},\mathbf{y}^{(t+1)},\mathbf{a}) \ \ \text{s.t.} \ \ \mathbf{a} \in [0,1]^n$. \; \\
        Update the regularization coefficient: $\lambda \leftarrow \lambda \times \rho$. \;\\
        Set $t \leftarrow t+1$.
        }
        \caption{Alternating minimization for biconvex case}        
        \label{alg:biconvex}
    \end{algorithm}
    
    \noindent \textbf{Remark (Extension to block multi-convex functions).} 
    \textit{Algorithm \ref{alg:biconvex} can be extended to the block multi-convex functions \cite{shen2017disciplined}. In this case, we update the arguments one by one at each iteration. Moreover, if the objective function is convex with respect to a subset of its arguments, we can update the whole subset at once.} 
    
    In the following lemma, we prove that if $f(\mathbf{x},\mathbf{y})$ is $L$-\textit{Lipschitz} on $\mathbf{x}$, then $L > \lambda^\star$. As a result, $\lambda > L$ guarantees $\lambda > \lambda^\star$.
    
     \noindent \textbf{Lemma 2.} 
    \textit{Assume that $f(\mathbf{x},\mathbf{y})$ is an L-Lipschitz function on $\mathbf{x} \in [0,1]^n$. Then, $\lambda > L$ guarantees $\phi_{\lambda}(\mathbf{x}^{\star},\mathbf{a}^{\star}) = 0$, where}
    \begin{align*}
        (\mathbf{x}^{\star},\mathbf{a}^{\star}) = \underset{\mathbf{x},\mathbf{a}}{\text{argmin}} \ \mathcal{L}_{\lambda} (\mathbf{x},\mathbf{y}^{(t)},\mathbf{a}) \ \ 
        \text{s.t.} \ \ \mathbf{x},\mathbf{a}\in [0,1]^n , \ \mathbf{x}\in \mathcal{D}.
    \end{align*}
    \noindent \textit{Proof.} If $\mathbf{x}^\star \in\{0,1\}^n$, based on \eqref{eq:update-a}, $\mathbf{a}^\star=\widetilde{\mathbf{x}}^\star=\mathbf{x}^\star$ and then, $\phi_{\lambda}(\mathbf{x}^\star,\mathbf{a}^\star) = 0$. Otherwise, substituting $\mathbf{a}^\star=\widetilde{\mathbf{x}}^\star$ in $\mathcal{L}_{\lambda} (\mathbf{x},\mathbf{y}^{(t)},\mathbf{a})$ results in the following sub-problem for $\mathbf{x}$
    \begin{align*}
        \mathbf{x}^\star = \ \underset{\mathbf{x}}{\text{argmin}} f(\mathbf{x},\mathbf{y}^{(t)}) + \lambda d(\mathbf{x}) \ \
        \text{s.t.} \ \ \mathbf{x}\in [0,1]^n \cap \mathcal{D}. 
    \end{align*}
    By dividing the elements of vector $\mathbf{x}$ into two groups, $x_i \geq 1/2$ and $x_i < 1/2$, $\forall i$, it can be shown that $\| \widetilde{\mathbf{x}}-\mathbf{x}\| \leq d(\mathbf{x}), \ \ \forall \mathbf{x} \in [0,1]^n$. Now, we have
    \begin{align}
    \label{eq:lemma3}
       \lambda \| \widetilde{\mathbf{x}}^\star-\mathbf{x}^\star\| \leq \  \lambda  d(\mathbf{x}^\star) \nonumber = \ & \lambda d(\mathbf{x}^\star) + f(\mathbf{x}^\star,\mathbf{y}^{(t)}) - f(\mathbf{x}^\star,\mathbf{y}^{(t)})\nonumber  \\
       = \ & \underset{\mathbf{x}}{\min} \left\{f(\mathbf{x},\mathbf{y}^{(t)}) + \lambda d(\mathbf{x})\right\} - f(\mathbf{x}^\star,\mathbf{y}^{(t)}) \nonumber \\
       \leq \ &\lambda d(\widetilde{\mathbf{x}}^\star) + f(\widetilde{\mathbf{x}}^\star,\mathbf{y}^{(t)}) - f(\mathbf{x}^\star,\mathbf{y}^{(t)}) \nonumber  \\
       = \ & f(\widetilde{\mathbf{x}}^\star,\mathbf{y}^{(t)}) - f(\mathbf{x}^\star,\mathbf{y}^{(t)}) \leq \ L \| \widetilde{\mathbf{x}}^\star-\mathbf{x}^\star\|.
    \end{align}
    \noindent From \eqref{eq:lemma3}, we have $\| \widetilde{\mathbf{x}}^\star-\mathbf{x}^\star\| (\lambda - L) \leq 0$. Thus, based on the assumption $\lambda > L$, we conclude that $\| \widetilde{\mathbf{x}}^\star-\mathbf{x}^\star\| = 0$, meaning $\mathbf{x}^\star = \widetilde{\mathbf{x}}^\star  \in \{0,1\}^n$ and hence, $\phi_{\lambda}(\mathbf{x}^{\star},\mathbf{a}^{\star}) = 0$.
    \hfill\(\blacksquare\) \\
    
    \noindent Lemma 2 shows that by choosing $\lambda > L$, the solution of Algorithm \ref{alg:biconvex} for variable $\mathbf{x}$ is binary. Now, in the following theorem, we express the convergence rate of Algorithm \ref{alg:biconvex}. 
    
    \noindent \textbf{Theorem 2 (Convergence rate).} 
    \textit{Assume that $f(\mathbf{x},\mathbf{y})$ is an L-Lipschitz function on $\mathbf{x} \in [0,1]^n$. For $\lambda \geq L\sqrt{n} / \epsilon$, Algorithm \ref{alg:biconvex} converges to $(\mathbf{x}^\star,\mathbf{y}^\star,\mathbf{a}^\star)$ with $d(\mathbf{x}^\star) \leq \epsilon$ in at most $\ceil{\left(\log (L\sqrt{n}) - \log (\epsilon \lambda^{(0)}) \right) /  \log \rho}$ iterations where $\lambda^{(0)}$ is the initial value of $\lambda$ and $\ceil{.}$ is the ceiling operator.}
    
    \noindent \textit{Proof.} 
    We start with the variational inequality \cite{he20121} which states that for every $\mathbf{x}\in[0,1]^n \cap  \mathcal{D}$, we have
    \begin{align*}
    \frac{1}{\lambda}\langle\mathbf{x}-\mathbf{x}^\star, \nabla _{\mathbf{x}}f(\mathbf{x}^\star,\mathbf{y})\rangle + 
    d(\mathbf{x}) - d(\mathbf{x}^\star) \geq 0.
    \end{align*}
    \noindent For any feasible $\mathbf{x} \in \{ 0,1 \}^n \cap \mathcal{D}$, we have
    \begin{align}
    \label{eq:convergence}
     d(\mathbf{x}^\star) 
     &\leq \overbrace{d(\mathbf{x})}^{=0}  + \frac{1}{\lambda}\langle\mathbf{x}-\mathbf{x}^\star, \nabla_{\mathbf{x}}f(\mathbf{x}^\star,\mathbf{y})\rangle \nonumber \\
     &\leq \frac{1}{\lambda} \left\|\mathbf{x}-\mathbf{x}^\star \right\| \left \| \nabla_{\mathbf{x}}f(\mathbf{x}^\star,\mathbf{y}) \right\| \leq \frac{L\sqrt{n}}{\lambda} \leq \epsilon,
    \end{align}
    where the first inequality comes from the mentioned variational inequality, the second one holds due to Cauchy-Schwarz inequality, and the third one comes from the Lipschitz continuity assumption and the fact that $\left\|\mathbf{x}-\mathbf{x}^\star \right\| \leq \sqrt{n}$. Now, at the $t$-th iteration, the updated regularization coefficient is $\lambda^{(t)} = \rho^t\lambda^{(0)}$. Thus, we have
    $\lambda^{(t)} = \rho^t\lambda^{(0)} \geq L\sqrt{n} / \epsilon$, resulting in $t \geq \left(\log (L\sqrt{n}) - \log (\epsilon \lambda^{(0)}) \right) /  \log \rho$.
    \hfill\(\blacksquare\) 
\subsection{General case}
\textcolor{black}{Now, we extend our algorithm to the case in which $f(\mathbf{x},\mathbf{y})$ and feasible region $\mathcal{D}$ are not necessarily biconvex and convex, respectively.} In this case, the proposed algorithm can not be directly used. This is due to the fact that the corresponding sub-problems in Algorithm \ref{alg:biconvex} are not convex over $\mathbf{x}$ and $\mathbf{y}$. In this case, we need to modify problem \eqref{eq:deductive-penalty} prior to utilizing Algorithm \ref{alg:biconvex}. {According to \cite{tuy1998convex}, every twice continuously differentiable function can be written as the difference of two convex functions.} As a result, $f(\mathbf{x},\mathbf{y})$ can be expressed as
\begin{align*}
f(\mathbf{x},\mathbf{y}) = f_a(\mathbf{x},\mathbf{y}) - f_b(\mathbf{x},\mathbf{y}),    
\end{align*}
where $f_{a}(\mathbf{x},\mathbf{y})$ and $f_{b}(\mathbf{x},\mathbf{y})$ are two convex functions with respect to $\mathbf{x}$ and $\mathbf{y}$. We use this fact to develop an iterative method to solve \eqref{eq:deductive-penalty}. Let $l$ denote the iteration index. In the $(l+1)$-th iteration of this method,
$f_{b}(\mathbf{x},\mathbf{y})$ is linearized with 
\begin{align}
\label{f_B_bar:DC}
&{\bar{f}}_{b}^{(l+1)}(\mathbf{x},\mathbf{y}; \mathbf{x}^{(l)},\mathbf{y}^{(l)}) := f_{b}(\mathbf{x}^{(l)},\mathbf{y}^{(l)}) + \left\langle \nabla_{\mathbf{x}} f_{b} (\mathbf{x}^{(l)},\mathbf{y}^{(l)}), \mathbf{x} - \mathbf{x}^{(l)} \right\rangle +  \left\langle\nabla_{\mathbf{y}} f_{b} (\mathbf{x}^{(l)},\mathbf{y}^{(l)}), \mathbf{y} - \mathbf{y}^{(l)}\right\rangle,
\end{align}
where $\mathbf{x}^{(l)}$ and $\mathbf{y}^{(l)}$ are the solutions of the $l$-th iteration. Using \eqref{f_B_bar:DC}, in the $(l+1)$-th iteration, the objective function is replaced with 
\begin{equation*}
\bar{f}^{(l+1)}(\mathbf{x},\mathbf{y}; \mathbf{x}^{(l)},\mathbf{y}^{(l)}) : =f_{a}(\mathbf{x},\mathbf{y}) - {\bar{f}}_{b}^{(l+1)}(\mathbf{x},\mathbf{y}; \mathbf{x}^{(l)},\mathbf{y}^{(l)}),
\end{equation*}
\noindent \textcolor{black}{which is a convex function in $\mathbf{x}$ and $\mathbf{y}$.} In the case that $\mathcal{D}$ is not convex, we can adopt a similar approach for the non-convex constraints $g_i(\mathbf{x},\mathbf{y})$. In other words, we can write 
\begin{align*}
g_i(\mathbf{x},\mathbf{y}) = g_{i,a}(\mathbf{x},\mathbf{y}) - g_{i,b}(\mathbf{x},\mathbf{y}),    
\end{align*}
where $g_{i,a}(\mathbf{x},\mathbf{y})$ and $g_{i,b}(\mathbf{x},\mathbf{y})$ are convex functions with respect to $\mathbf{x}$ and $\mathbf{y}$. In the $(l+1)$-th iteration, constraint $g_i(\mathbf{x},\mathbf{y}) \leq 0$ is adjusted to 
\begin{align*}
g_{i,a}(\mathbf{x},\mathbf{y}) - {\bar{g}}_{i,b}^{(l+1)}(\mathbf{x},\mathbf{y}; \mathbf{x}^{(l)},\mathbf{y}^{(l)}) \leq 0,    
\end{align*}
where 
${\bar{g}}_{i,b}^{(l+1)}(\mathbf{x},\mathbf{y}; \mathbf{x}^{(l)},\mathbf{y}^{(l)})$ is the linearized version of ${g}_{i,b}(\mathbf{x},\mathbf{y})$ derived similar to \eqref{f_B_bar:DC}. If $\bar{\mathcal{D}}^{(l+1)}$ denotes the convexified feasible region of the problem in the $(l+1)$-th iteration, the problem is expressed as
\begin{align}
	\label{eq:DC:original2}
	&\underset{\mathbf{x},\mathbf{y}}{\min} 
	\ \ \ \ \ \ \bar{f}^{(l+1)}(\mathbf{x},\mathbf{y}; \mathbf{x}^{(l)},\mathbf{y}^{(l)})   \\
	&\text{s.t.}  
	\ \ \ \ \ \ \ \mathbf{x},\mathbf{y} \in \bar{\mathcal{D}}^{(l+1)}, \ \ \mathbf{x} \in \{0 , 1\}^n. \nonumber
\end{align}
We can solve this problem via Algorithm \ref{alg:biconvex}. To solve \eqref{eq:deductive-penalty}, we start at $l=0$ with a feasible point $(\mathbf{x}^{(l)},\mathbf{y}^{(l)})$, and iteratively solve \eqref{eq:DC:original2} via Algorithm \ref{alg:biconvex}. The solution obtained by Algorithm \ref{alg:biconvex} is used to update $\bar{f}^{(l+1)}(\mathbf{x},\mathbf{y}; \mathbf{x}^{(l)},\mathbf{y}^{(l)})$ and $\bar{\mathcal{D}}^{(l+1)}$
 in the next iteration. This procedure is presented in Algorithm \ref{alg:nonconvex}. It can be shown that Algorithm \ref{alg:nonconvex} generates a sequence of improved solutions \cite{BehradCluster}, i.e.,
 \begin{align*}
 f(\mathbf{x}^{(l)},\mathbf{y}^{(l)}) +  \phi_{\lambda}(\mathbf{x}^{(l)},\mathbf{a}^{(l)}) \geq     f(\mathbf{x}^{(l+1)},\mathbf{y}^{(l+1)}) +  \phi_{\lambda}(\mathbf{x}^{(l+1)},\mathbf{a}^{(l+1)}).    
 \end{align*}
 {Moreover, the solution of problem \eqref{eq:DC:original2} is feasible for the main problem as $\bar{\mathcal{D}}^{(l+1)} \subseteq \mathcal{D}, \forall l$ \textcolor{black}{(see Appendix \ref{app:trust})}.} \\
 
 \noindent Using Theorem $2$, we can obtain the convergence rate of Algorithm \ref{alg:nonconvex} when $f_{a}(\mathbf{x},\mathbf{y})$ is an $L$-\textit{Lipschitz} function on $\mathbf{x} \in [0,1]^n$. We first explain the following useful lemma.

    \noindent \textbf{Lemma 3.} 
    \textit{Consider the objective function in \eqref{eq:DC:original2}. If $f_{a}(\mathbf{x},\mathbf{y})$ is an $L$-Lipschitz function on $\mathbf{x} \in [0,1]^n$, $\bar{f}^{(l+1)}(\mathbf{x},\mathbf{y}; \mathbf{x}^{(l)},\mathbf{y}^{(l)})$ also satisfies Lipschitz continuity on $\mathbf{x} \in [0,1]^n$ with constant $L + \left\| \nabla_{\mathbf{x}} f_{b} (\mathbf{x}^{(l)},\mathbf{y}^{(l)})\right\|$.}
    
    \noindent \textit{Proof.} \textcolor{black}{See Appendix \ref{app:LipDC}.}    \hfill\(\blacksquare\).
    
    \noindent Using lemma 3 and Theorem 2, the number of inner iterations required for convergence of Algorithm \ref{alg:nonconvex} at the $(l+1)$-th iteration is $\ceil{\left(\log (L^\prime\sqrt{n}) - \log (\epsilon \lambda^{(0)}) \right) /  \log \rho)}$
    where $L^\prime = L + \left\| \nabla_{\mathbf{x}} f_{b} (\mathbf{x}^{(l)},\mathbf{y}^{(l)})\right\|$.

   \begin{algorithm}[H]
        \SetAlgoLined
        Set $l = 0$ and initialize $\mathbf{x}^{(l)}$, $\mathbf{y}^{(l)}$.\;\\
        \Repeat{convergence of $\mathbf{x}$ and $\mathbf{y}$}{
        Construct $\bar{\mathcal{D}}^{(l+1)}$. \; \\
        Solve problem \eqref{eq:DC:original2} via Algorithm \ref{alg:biconvex}.\; \\
        Update $(\mathbf{x}^{(l+1)}$, $\mathbf{y}^{(l+1)})$. \; \\
        Set $l \leftarrow l+1$. \; 
        }
        \caption{Alternating minimization for general case}
        \label{alg:nonconvex}
    \end{algorithm}
    
    \section{Simulation Studies}

    In this section, we apply our algorithm to a RAT selection problem in HetNets. We propose a user-RAT assignment strategy in a multi-RAT network equipped by two technologies of WiFi (RAT-1) and OFDMA (RAT-2). Our goal is to maximize the network aggregate throughput \cite{keshavarz2016hetnets}. 
    
    \noindent {\it{Problem definition}}: We consider a HetNet including $I$ users and $K$ small base stations (SBSs). The set of users and SBSs are presented by $\mathcal{I}$ and $\mathcal{K}$, respectively. We use indices $i$ and $k$ to denote users and SBSs. Each SBS belongs to one of two different RATs: \textit{WiFi} (RAT-1) and \textit{OFDMA-based cellular} (RAT-2). The set of RAT-$m$ SBSs and the set of users connected to RAT-$m$ are illustrated by $\mathcal{K}_m$ and $\mathcal{I}_m$, $m=1,2$, respectively. Each user connects to only one SBS. We denote the assignment of the $i$-th user to the $k$-th SBS by the binary variable $x_{ik}$ being 1. Otherwise, $x_{ik}=0$. Using this notation, the $i$-th user connects to RAT-$m$, if and only if $\sum_{k \in \mathcal{K}_m}x_{ik} = 1, \ m=1,2$. We present the instantaneous rate of the $i$-th user connected to the $k$-th SBS by $r_{ik}$.
 
    \noindent Let $w_{i}^{(m)}$ denote the average throughput of the $i$-th user when connected to RAT-$m$. According to \cite{keshavarz2016hetnets}, we can write 
    \begin{align*}
     &w^{(1)}_i = \left( \sum_{j\in\mathcal{I}}  \sum_{k \in \mathcal{K}_1} \frac{x_{jk}}{r_{jk}} \right)^{-1}, \ \forall i \in \mathcal{I}_1, \\ 
     &w^{(2)}_i = \left({\sum\limits_{k \in \mathcal{K}_2} x_{ik}r_{ik}}\right) {\left(\sum\limits_{j \in \mathcal{I}} \sum\limits_{k \in \mathcal{K}_2}x_{jk}\right)^{-1}}, \ \forall i \in \mathcal{I}_2.
    \end{align*}
    Thus, the throughput of the $i$-th user is 
    \begin{align*}
    w_i = \sum_{k \in \mathcal{K}_1}x_{ik} w^{(1)}_i +\sum_{k \in \mathcal{K}_2}x_{ik}  w^{(2)}_i.
    \end{align*}
    \noindent The corresponding optimization problem to maximize the aggregate throughput can be expressed as
   \begin{align}
	\label{RAT-problem}
	&\underset{\mathbf{x}}{\max} \hspace{2mm}
	\ \ \ \ \ \ \sum\limits_{i \in \mathcal{I}} \alpha_i w_i \\
	&\text{s.t.}  
	\ \ \ \ \ \ \ \ \ \ \text{C1:}  \ \ \sum\limits_{i \in \mathcal{I}} \sum\limits\limits_{k \in \mathcal{K}_m}x_{ik} w^{(m)}_i \leq w^{(m)}_{\max},  \ \ m=1,2,\nonumber \\
	& \hspace{2cm} \text{C2:}  \ \ \sum\limits_{i \in \mathcal{I}}\sum\limits_{k \in \mathcal{K}_m} x_{ik} \leq N^{(m)}_{\max},\ \ m=1,2,\nonumber \\
	&\hspace{2cm}\text{C3:} \ \ \sum\limits_{k \in \mathcal{K}} x_{ik} = 1, \ \ \forall i \in \mathcal{I},\nonumber \\
	& \hspace{2cm}\text{C4:} \ \ \sum\limits_{i \in \mathcal{I}} x_{ik} \leq K_{\max}^{(k)}, \ \ \forall k \in \mathcal{K},\nonumber \\
	& \hspace{2cm}\text{C5:} \ \ x_{ik} \in \{0,1\}, \ \ \forall k \in \mathcal{K},\forall i \in \mathcal{I},\nonumber 
	\end{align}
    where $\alpha_i>0, \ \forall i$ is the fairness coefficient between users. Moreover, $w^{(m)}_{\max}$ and $N^{(m)}_{\max}$ present the maximum throughout and the maximum number of users connected to RAT-$m$, respectively, and $K_{\max}^{(k)}$ is the maximum number of users served by the $k$-th SBS. In \eqref{RAT-problem}, C1 and C2 together imply that the throughput of the RATs is limited to the backhaul capacity. In particular, C1 restricts the throughput of RATs, and C2 confines the number of users connected to the RATs. Constraint C3 guarantees that each user connects to exactly one SBS. Constraint C4 shows that the number of users connected to each SBS is limited. Finally, C5 denotes the association indicators are binary.
    
    \textcolor{black}{In what follows, we propose two solutions for problem \eqref{RAT-problem} based on the proposed algorithms.}
    \begin{itemize}
    \item \textbf{Algorithm \ref{alg:biconvex}:}
    To use Algorithm \ref{alg:biconvex}, the objective function must be biconvex (block multi-convex). However, the objective function in \eqref{RAT-problem} does not satisfy this requirement. To tackle this issue, we use the following change of variables 
    \begin{align}
    \label{eq:new-var}
    v_1 = \sum\limits\limits_{i\in\mathcal{I}} \sum\limits\limits_{k \in \mathcal{K}_1} \frac{x_{ik}}{r_{ik}},
	v_2=\sum\limits\limits_{i \in \mathcal{I}} \sum\limits\limits_{k \in \mathcal{K}_2} x_{ik},
	u_i = \sum\limits_{k\in \mathcal{K}_2} x_{ik}r_{ik}.
    \end{align}
    \normalsize{}
    \noindent Let $\mathbf{u}=[u_i, \forall i]$ and $\mathbf{v}=[v_1,v_2]$. The modified objective function is 
    \begin{align*}
        f(\mathbf{x},\mathbf{u},\mathbf{v})=\sum\limits_{i \in \mathcal{I}} \frac{\alpha_i}{v_1} \sum\limits_{k\in \mathcal{K}_1}x_{ik} + \sum\limits_{i \in \mathcal{I}} \frac{\alpha_i}{v_2} \sum\limits_{k\in \mathcal{K}_2}x_{ik}u_i,
    \end{align*}
    and the corresponding optimization problem is given by
    \begin{align}
	\label{RAT-bivoncex}
	\underset{\mathbf{x},\mathbf{u},\mathbf{v}}{\max} \ \ \
	 f(\mathbf{x},\mathbf{u},\mathbf{v}) 
	 \ \ \ \text{s.t.} \hspace{4mm}
	 \text{C1-C5} , \eqref{eq:new-var},\ \mathbf{u},\mathbf{v} \geq \mathbf{0}.
	\end{align}
    It can be shown that $f(\mathbf{x},\mathbf{u},\mathbf{v})$ is a block multi-convex function with respect to its arguments. Thus, we can employ Algorithm \ref{alg:biconvex} to solve \eqref{RAT-bivoncex}.
	Moreover, using the Cauchy-Schwarz inequality, one can show that $f(\mathbf{x},\mathbf{u},\mathbf{v})$ is an $L$-Lipschitz function on $\mathbf{x} \in [0,1]^{I\times K}$ with 
	\begin{align*}
	&\hspace{12mm}L=\sqrt{I}\alpha_{\max}\left| {r^{(1)}_{\max}}-{r^{(2)}_{\min}}{(I-N^{(1)}_{\max})}^{-1} \right|, \\
	&\alpha_{\max} = \underset{i \in \mathcal{I}}{\max} \ \alpha_i, r^{(1)}_{\max} = \underset{i \in \mathcal{I},k \in \mathcal{K}_1}{\max} \ r_{ik},r^{(2)}_{\min} = \underset{i \in \mathcal{I},k \in \mathcal{K}_2}{\min} \ r_{ik}.
	\end{align*}
	
    \item \textbf{Algorithm \ref{alg:nonconvex}:}
    In order to use Algorithm \ref{alg:nonconvex}, we just need to write the objective function of problem \eqref{RAT-problem} as the difference of two convex functions. Using the change of variables given in \eqref{eq:new-var}, problem \eqref{RAT-problem} is equivalent to \eqref{RAT-bivoncex}. In the following lemma, we show that $f(\mathbf{x},\mathbf{u},\mathbf{v})$ can be written as the difference of two convex functions.
    
    \textbf{Lemma 4.} \textit{Function $f(a,b,c)=\frac{ab}{c}=\frac{(a+b)^2}{2c} - (\frac{a^2}{2c} + \frac{b^2}{2c})$ is written as the difference of two convex functions for $a,b,c >0$.}
    
    \textit{Proof.} According to \cite{boyd2004convex}, chapter 3, if $f(\mathbf{x})$ is convex, then $g(\mathbf{x},t)=tf(\mathbf{x}/t)$ is also convex on $\{(\mathbf{x},t)\big| \mathbf{x}/t \in \text{dom} \ f\}$. Since $(a+b)^2$, $a^2$ , and $b^2$ are convex functions for $a,b,c >0$, functions $\frac{(a+b)^2}{2c}$, $\frac{a^2}{2c}$, and $\frac{b^2}{2c}$ are also convex. As a results, $f(a,b,c)=\frac{ab}{c}$ can be written as difference of two convex functions. 
    \hfill\(\blacksquare\) \\
    Using the above lemma, we can write $f(\mathbf{x},\mathbf{u},\mathbf{v})=f_a(\mathbf{x},\mathbf{u},\mathbf{v}) - f_b(\mathbf{x},\mathbf{u},\mathbf{v})$, where
    
    \begin{align*}
        f_a(\mathbf{x},\mathbf{u},\mathbf{v})=\frac{1}{2}\sum\limits_{i \in \mathcal{I}} {\alpha_i} \left[ \sum\limits_{k\in \mathcal{K}_1}\left( x_{ik} + \frac{1}{v_1}\right)^2 + \sum\limits_{k\in \mathcal{K}_2} \frac{\left( x_{ik} + u_i \right)^2 }{v_2} \right],    
    \end{align*}
    and
    \begin{align*}
        f_b(\mathbf{x},\mathbf{u},\mathbf{v})=\frac{1}{2}\sum\limits_{i \in \mathcal{I}} {\alpha_i} \left[ \sum\limits_{k\in \mathcal{K}_1}\left( (x_{ik})^2 + \frac{1}{(v_1)^2}\right) + \sum\limits_{k\in \mathcal{K}_2} \frac{( x_{ik})^2 + (u_i)^2 }{v_2} \right],    
    \end{align*}
are two convex functions. It is worth mentioning that we use the fact that sum of convex functions is also convex and $v_1,v_2 > 0$.

We can write the linearized version of $f(\mathbf{x},\mathbf{u},\mathbf{v})$ at the $l$-th iteration as
\begin{equation*}
\bar{f}^{(l+1)}(\mathbf{x},\mathbf{u},\mathbf{v}; \mathbf{x}^{(l)},\mathbf{u}^{(l)},\mathbf{v}^{(l)}) =f_{a}(\mathbf{x},\mathbf{u},\mathbf{v}) - {\bar{f}}_{b}^{(l+1)}(\mathbf{x},\mathbf{u},\mathbf{v}; \mathbf{x}^{(l)},\mathbf{u}^{(l)},\mathbf{v}^{(l)}),
\end{equation*}
 where
\begin{align*}
{\bar{f}}_{b}^{(l+1)}(\mathbf{x},\mathbf{u},\mathbf{v}; \mathbf{x}^{(l)},\mathbf{u}^{(l)},\mathbf{v}^{(l)}) &= f_{b}(\mathbf{x}^{(l)},\mathbf{u}^{(l)},\mathbf{v}^{(l)})  + \left\langle \nabla_{\mathbf{x}} f_{b} (\mathbf{x}^{(l)},\mathbf{u}^{(l)},\mathbf{v}^{(l)}), \mathbf{x} - \mathbf{x}^{(l)} \right\rangle \\
& +  \left\langle\nabla_{\mathbf{u}} f_{b} (\mathbf{x}^{(l)},\mathbf{u}^{(l)},\mathbf{v}^{(l)}), \mathbf{u} - \mathbf{u}^{(l)}\right\rangle + \left\langle\nabla_{\mathbf{v}} f_{b} (\mathbf{x}^{(l)},\mathbf{u}^{(l)},\mathbf{v}^{(l)}), \mathbf{v} - \mathbf{v}^{(l)}\right\rangle. 
\end{align*}
    
    \end{itemize}
    
	
	\begin{figure}
        \centering
        \includegraphics[trim={0cm 0 0cm 0cm},width=0.6\linewidth]{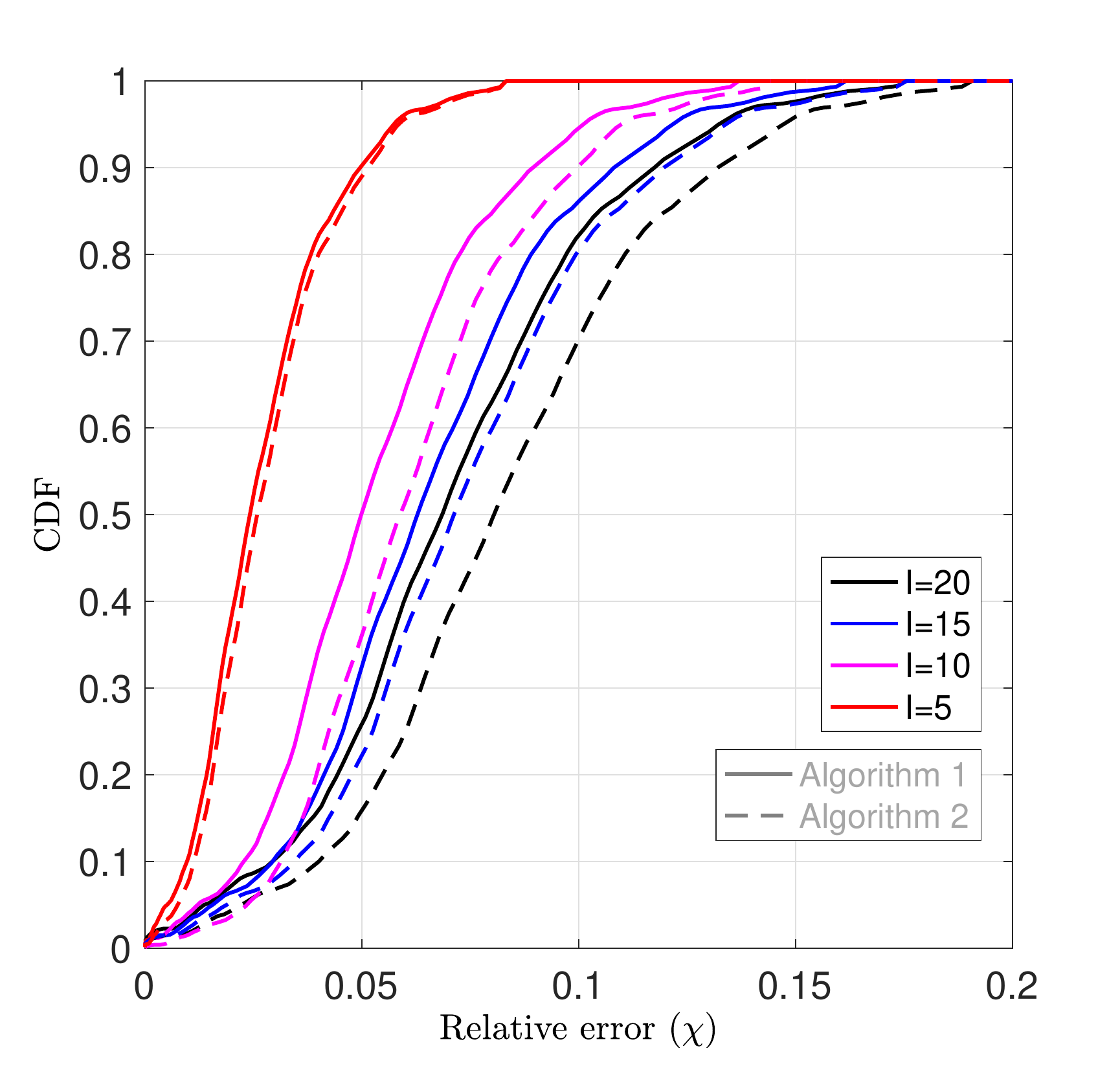}
        \caption{CDF of relative error with optimal solution.}
        \label{fig:CDF}
    \end{figure}

	\noindent {\it{Performance evaluation}}:
	We consider a cell with radius $50$m including $K=2$ SBSs. Each SBS belongs to a different technology (RAT-1 or -2). The bandwidth and the noise power density are $180$kHz and $-174$dBm/Hz, respectively. The channel gain between the $i$-th user and the $k$-th SBS is expressed as $h_{ik}=\psi d_{ik}^{-3}$ where $d_{ik}$ is the distance between the $i$-th user and the $k$-th SBS and $\psi$ is a Rayleigh random variable with unit variance. The parameters of problem \eqref{RAT-problem} are $	N_{\max}^{(m)}=I-1,  w_{\max}^{(m)}=\underset{i \in \mathcal{I},k \in \mathcal{K}_m}{\max} \ r_{ik},  K^{(k)}_{\max}=I (\forall k), \alpha_i=1 (\forall i)$.
	
	To find the sub-optimality of the proposed algorithms, we calculate the relative error as $\chi = 1-\frac{f^\star}{f_\text{alg.}}$ where $f^\star$ is the optimal value of the objective function corresponding to the exhaustive search, and $f_\text{alg.}$ is the solution of the algorithms. Fig. \ref{fig:CDF} shows the cumulative distribution function (CDF) of $\chi$ for $100$ different realizations corresponding to $I = 5, 10, 15$, and $20$. As presented, $85\%$ sub-optimality has been achieved with a probability more than $95\%$.
	
	Fig. \ref{fig:lambda} depicts the solution of the regularized problem in \eqref{RAT-bivoncex} versus the regularization coefficient $\lambda$. We assume $I = 10,50$, and $100$, and $\rho = 1$. The Lipschitz constants are shown with dashed lines. As proven in Lemma 2, $\lambda > L$ ensures the solution is binary. After this point, increasing the value of $\lambda$ does not improve the value of $f_{\text{alg.}}$. {This demonstrates the inefficiency of cross-validation techniques used for MIPs.} 

    Fig. \ref{fig:rate} illustrates the convergence of Algorithm \ref{alg:biconvex} for different values of $\rho$ and $\epsilon$. We assume $\lambda^{(0)} = 1$ and $I = 50$. The vertical lines are the theoretical convergence bounds derived in Theorem 2. 
    As can be seen, all curves converge prior to their corresponding theoretical bounds. 
	
		\begin{figure}
        \centering
        \includegraphics[trim={0cm 0 0cm 0cm},width=0.6\linewidth]{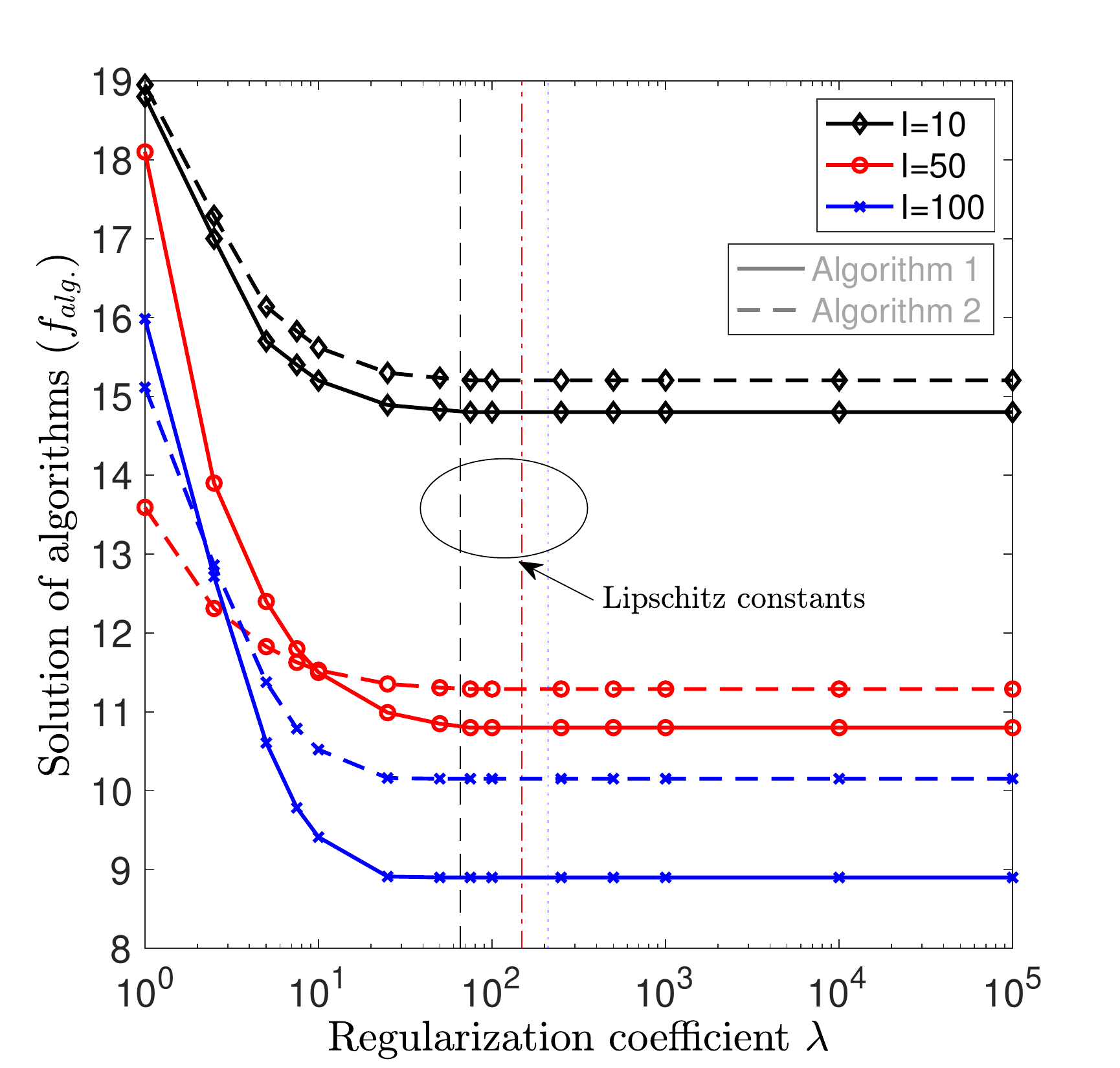}
        \caption{The value of $f_{\text{alg.}}$ versus $\lambda$ (the curves with the same colors correspond to the same scenario).}
        \label{fig:lambda}
    \end{figure}
    
    \begin{figure}
        \centering
        \includegraphics[trim={0cm 0 0cm 0cm},width=0.6\linewidth]{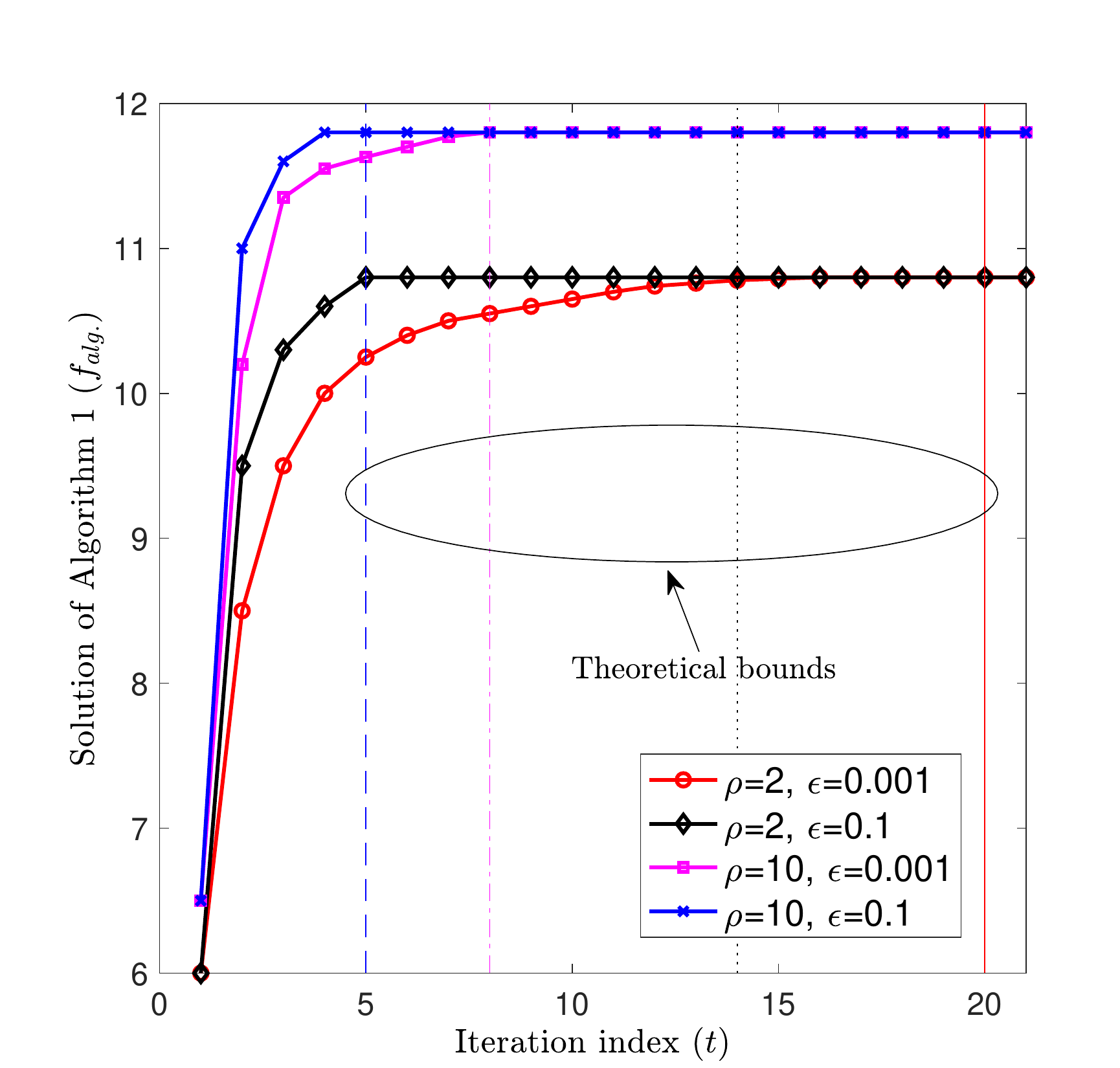}
        \caption{Convergence of Algorithm \ref{alg:biconvex} (the curves with the same colors correspond to the same scenario).}
        \label{fig:rate}
    \end{figure}

	\section{Conclusion}
	In this paper, we proposed a regularized alternating method to solve MIPs. Instead of using cross-validation to find a proper regularization coefficient, we update its value iteratively. We have shown that our algorithm results in a sequence of improved solutions. Moreover, for the case that the objective function is continuously differentiable, we derived the convergence rate, a lower bound on the value of regularization coefficient, and an upper bound on the number of iterations required for the convergence.

    \appendices
    \section{$\omega(\lambda)$ is increasing}\label{app:omega}	
    To show that function $\omega(\lambda)$ is increasing with respect to $\lambda$, first we need to show that $\mathcal{L}_{\lambda} (\mathbf{x},\mathbf{y},\mathbf{a})$ is an increasing function with respect to $\lambda$. According to the Cauchy-Schwarz inequality (Lemma 1), 
\begin{align*}
\left(n - \langle\mathbf{x},\mathbf{a}\rangle - \langle\mathbf{1}-\mathbf{x},\mathbf{1}-\mathbf{a}\rangle\right) \geq 0, \ \ \forall \mathbf{x},\mathbf{a}\in [0,1]^n.
\end{align*}
Let $\lambda_1 > \lambda_2 > 0$. Since
\begin{align*}
    \mathcal{L}_{\lambda} (\mathbf{x},\mathbf{y},\mathbf{a}) = f(\mathbf{x},\mathbf{y}) + \lambda \left(n - \langle\mathbf{x},\mathbf{a}\rangle - \langle\mathbf{1}-\mathbf{x},\mathbf{1}-\mathbf{a}\rangle\right),
\end{align*}
we can write 
\begin{align*}
    \mathcal{L}_{\lambda_1} (\mathbf{x},\mathbf{y},\mathbf{a})  \geq \mathcal{L}_{\lambda_2} (\mathbf{x},\mathbf{y},\mathbf{a}), \hspace{1cm} \forall (\mathbf{x},\mathbf{y},\mathbf{a}) \in \mathcal{R}.
\end{align*}
As a result, function $\mathcal{L}_{\lambda} (\mathbf{x},\mathbf{y},\mathbf{a})$ is increasing with respect to $\lambda$. According to the definition of $\omega(\lambda)$, we have 
\begin{align*}
    \omega(\lambda) := \underset{(\mathbf{x},\mathbf{y},\mathbf{a})\in \mathcal{R}}{\min} \ \mathcal{L}_{\lambda} (\mathbf{x},\mathbf{y},\mathbf{a}).
\end{align*}
Hence, we can write
\begin{align*}
    \omega(\lambda_1) = \underset{(\mathbf{x},\mathbf{y},\mathbf{a})\in \mathcal{R}}{\min} \ \mathcal{L}_{\lambda_1} (\mathbf{x},\mathbf{y},\mathbf{a}), \\
    \omega(\lambda_2) = \underset{(\mathbf{x},\mathbf{y},\mathbf{a})\in \mathcal{R}}{\min} \ \mathcal{L}_{\lambda_2} (\mathbf{x},\mathbf{y},\mathbf{a}).
\end{align*}
Since for any $(\mathbf{x},\mathbf{y},\mathbf{a}) \in \mathcal{R}$, we have $\mathcal{L}_{\lambda_1} (\mathbf{x},\mathbf{y},\mathbf{a})  \geq \mathcal{L}_{\lambda_2} (\mathbf{x},\mathbf{y},\mathbf{a})$,  We can apply the minimization operator to both sides as 
\begin{align*}
    \underset{(\mathbf{x},\mathbf{y},\mathbf{a})\in \mathcal{R}}{\min} \ \mathcal{L}_{\lambda_1} (\mathbf{x},\mathbf{y},\mathbf{a}) \geq \underset{(\mathbf{x},\mathbf{y},\mathbf{a})\in \mathcal{R}}{\min} \ \mathcal{L}_{\lambda_2} (\mathbf{x},\mathbf{y},\mathbf{a}).
\end{align*}
This results in 
\begin{align*}
    \omega(\lambda_1) \geq \omega(\lambda_2).
\end{align*}
Therefore, function $\omega(\lambda)$ is an increasing function of $\lambda$.

    \section{Trust region}\label{app:trust}	
    In what follows, we show that the feasible region of problem \eqref{eq:DC:original2} is a subset of the feasible region of the main problem, i.e., $\bar{\mathcal{D}}^{(l+1)} \subseteq \mathcal{D}, \forall l$. As a result, the solution of problem \eqref{eq:DC:original2}
 is feasible for the main problem in \eqref{eq:raw-problem}. 
 
\noindent The feasible region of the main problem is given in \eqref{feasiblereg}. Without loss of generality, we assume $p_{\text{nc}}$ number of the inequality constraints are non-convex and the remaining $p-p_{\text{nc}}$ constraints are all convex. Let $\mathcal{I}_{\text{nc}}$ and $\mathcal{I}_{\text{c}}$ denote the set of indices corresponding to the non-convex and convex constraints, respectively. We have $|\mathcal{I}_{\text{nc}}|=p_{\text{nc}}$ and $|\mathcal{I}_{\text{c}}|=p - p_{\text{nc}}$. For each $i \in \mathcal{I}_{\text{nc}}$, we can write $g_i(\mathbf{x},\mathbf{y}) = g_{i,a}(\mathbf{x},\mathbf{y}) - g_{i,b}(\mathbf{x},\mathbf{y})$, where $g_{i,a}(\mathbf{x},\mathbf{y})$ and $g_{i,b}(\mathbf{x},\mathbf{y})$ are convex functions with respect to $\mathbf{x}$ and $\mathbf{y}$. At the $(l+1)$-th iteration, the feasible region $\bar{\mathcal{D}}^{(l+1)}$ can be expressed as
\begin{align*}
    \bar{\mathcal{D}}^{(l+1)} =  \Big\{(\mathbf{x},\mathbf{y}) \Big| & g_i(\mathbf{x},\mathbf{y}) \leq 0, \  i \in \mathcal{I}_{\text{c}}, \\  & g_{i,a}(\mathbf{x},\mathbf{y}) - {\bar{g}}_{i,b}^{(l+1)}(\mathbf{x},\mathbf{y}; \mathbf{x}^{(l)},\mathbf{y}^{(l)}) \leq 0, \ i \in \mathcal{I}_{\text{nc}}, \\
	    &h_i(\mathbf{x},\mathbf{y}) = 0, \  i=1,2,\dots,q, \Big \},
\end{align*}
where ${\bar{g}}_{i,b}^{(l+1)}(\mathbf{x},\mathbf{y}; \mathbf{x}^{(l)},\mathbf{y}^{(l)})$ is the linearized version of ${g}^b_i(\mathbf{x},\mathbf{y})$. In the following lemma, we prove that $\bar{\mathcal{D}}^{(l+1)} \subseteq \mathcal{D}, \forall l$. 

\noindent \textbf{Lemma 5}. \textit{The feasible region of problem \eqref{eq:DC:original2}, i.e., $\bar{\mathcal{D}}^{(l+1)}$, is a subset of the feasible region of the main problem. In other words, $\bar{\mathcal{D}}^{(l+1)} \subseteq \mathcal{D}, \forall l$.}

\noindent \textit{Proof}. Let $(\mathbf{x},\mathbf{y}) \in \bar{\mathcal{D}}^{(l+1)}$. We have 
\begin{align*}
    & g_i(\mathbf{x},\mathbf{y}) \leq 0, \  i \in \mathcal{I}_{\text{c}}, \\  & g_{i,a}(\mathbf{x},\mathbf{y}) - {\bar{g}}_{i,b}^{(l+1)}(\mathbf{x},\mathbf{y}; \mathbf{x}^{(l)},\mathbf{y}^{(l)}) \leq 0, \ i \in \mathcal{I}_{\text{nc}}, \\
	    &h_i(\mathbf{x},\mathbf{y}) = 0, \  i=1,2,\dots,q.
\end{align*}
Since for each $i \in \mathcal{I}_{\text{nc}}$, $g_{i,b}(\mathbf{x},\mathbf{y})$ is a convex function, its linearized version is an under-estimator of $g_{i,b}(\mathbf{x},\mathbf{y})$. In other words, 
\begin{align*}
    {\bar{g}}_{i,b}^{(l+1)}(\mathbf{x},\mathbf{y}; \mathbf{x}^{(l)},\mathbf{y}^{(l)}) \leq g_{i,b}(\mathbf{x},\mathbf{y}), \ \ \  \forall (\mathbf{x},\mathbf{y}), \forall i \in \mathcal{I}_{\text{nc}}.
\end{align*}
As a result, 
\begin{align*}
    g_i(\mathbf{x},\mathbf{y}) & =  g_{i,a}(\mathbf{x},\mathbf{y}) - g_{i,b}(\mathbf{x},\mathbf{y}) \\
    & \leq g_{i,a}(\mathbf{x},\mathbf{y}) - {\bar{g}}_{i,b}^{(l+1)}(\mathbf{x},\mathbf{y}; \mathbf{x}^{(l)},\mathbf{y}^{(l)}) \leq 0, \ \ \ \forall i \in \mathcal{I}_{\text{nc}}.
\end{align*}
Hence, it is guaranteed to satisfy constraint $g_i(\mathbf{x},\mathbf{y}) \leq 0$, $\forall i \in \mathcal{I}_{\text{nc}}$. This means that $(\mathbf{x},\mathbf{y})$ satisfies all the following constraints as well
\begin{align*}
    & g_i(\mathbf{x},\mathbf{y}) \leq 0, \  i \in \mathcal{I}_{\text{c}}, \\  & g_i(\mathbf{x},\mathbf{y}) \leq 0, \  i \in \mathcal{I}_{\text{nc}}, \\
	    &h_i(\mathbf{x},\mathbf{y}) = 0, \  i=1,2,\dots,q.
\end{align*}
As a result, $(\mathbf{x},\mathbf{y}) \in \mathcal{D}$, and hence, $\bar{\mathcal{D}}^{(l+1)} \subseteq \mathcal{D}, \forall l$.
\hfill\(\blacksquare\) \\

\noindent The above lemma shows that any feasible solution of problem \eqref{eq:DC:original2} is also feasible for the main problem. As a result, we do not need to check the trust region.

    \section{Lipschitz continuity of linearized approximation}\label{app:LipDC}	
    From the definition of ${\bar{f}}_{b}^{(l+1)}(\mathbf{x},\mathbf{y}; \mathbf{x}^{(l)},\mathbf{y}^{(l)})$ in equation (15) of the manuscript, for any $\mathbf{x}_1,\mathbf{x}_2 \in [0,1]^n$ we have
   \begin{align}
   \label{sag}
        &\left|\bar{f}^{(l+1)}(\mathbf{x}_1,\mathbf{y}; \mathbf{x}^{(l)},\mathbf{y}^{(l)}) -\bar{f}^{(l+1)}(\mathbf{x}_2,\mathbf{y}; \mathbf{x}^{(l)},\mathbf{y}^{(l)})\right| = \nonumber \\
        &\left| f_{a}(\mathbf{x}_1,\mathbf{y}) -  {\bar{f}}_{b}^{(l+1)}(\mathbf{x}_1,\mathbf{y}; \mathbf{x}^{(l)},\mathbf{y}^{(l)}) -f_{a}(\mathbf{x}_2,\mathbf{y}) + {\bar{f}}_{b}^{(l+1)}(\mathbf{x}_2,\mathbf{y}; \mathbf{x}^{(l)},\mathbf{y}^{(l)}) \right| = \nonumber\\
        &\left|f_{a}(\mathbf{x}_1,\mathbf{y}) - f_{a}(\mathbf{x}_2,\mathbf{y}) + \left\langle \nabla_{\mathbf{x}} f_{b} (\mathbf{x}^{(l)},\mathbf{y}^{(l)}), \mathbf{x}_2 - \mathbf{x}_1 \right\rangle \right| \labelrel\leq{ineq1} \nonumber \\
        &\left|f_{a}(\mathbf{x}_1,\mathbf{y}) - f_{a}(\mathbf{x}_2,\mathbf{y}) | + |\left\langle \nabla_{\mathbf{x}} f_{b} (\mathbf{x}^{(l)},\mathbf{y}^{(l)}), \mathbf{x}_2 - \mathbf{x}_1 \right\rangle \right| \labelrel\leq{ineq2} \nonumber\\
        &\left|f_{a}(\mathbf{x}_1,\mathbf{y}) - f_{a}(\mathbf{x}_2,\mathbf{y}) \right| + \left\| \nabla_{\mathbf{x}} f_{b} (\mathbf{x}^{(l)},\mathbf{y}^{(l)})\right\| \left\|\mathbf{x}_2 - \mathbf{x}_1 \right\|,
    \end{align}
    where inequality (a) holds due to the triangle inequality and inequality (b) results from the Cauchy-Schwarz inequality. Now, using the fact that 
    $f_{a}(\mathbf{x},\mathbf{y})$ is $L$-Lipschitz on $\mathbf{x} \in [0,1]^n$ and following the definition of Lipschitz continuity along with inequality \eqref{sag}, one can conclude that the Lipschitz constant of $\bar{f}^{(l+1)}(\mathbf{x},\mathbf{y}; \mathbf{x}^{(l)},\mathbf{y}^{(l)})$ on $\mathbf{x} \in [0,1]^n$ is $L + \left\| \nabla_{\mathbf{x}} f_{b} (\mathbf{x}^{(l)},\mathbf{y}^{(l)})\right\|$.

    \section{Choosing arbitrary ($\rho > 1$, $\lambda > 0$)}\label{app:arbit}	
Our algorithm includes two parameters, $\rho > 1$ and $\lambda > 0$, which need initialization. In what follows, we show that for any arbitrary choice of these parameters, our algorithm works. Hence, these parameters do not require cross-validation.

According to our update in Algorithm 1, the value of $\lambda$ at the $t$-th iteration is $\lambda^{(t)}=\rho ^ t \lambda^{(0)}$. This update guarantees at some iteration $t=T$, we have $\lambda^{(T)}=\rho ^ T \lambda^{(0)} > \lambda^\star$. In other words, the sufficient condition of Theorem 1 holds. As a result, the solution of $\mathbf{x}$ becomes binary. In the following lemma, we show that for any arbitrary choice of ($\rho > 1$,$\lambda^{(0)} >0$), there is a $t=T$ such that $\lambda^{(T)} > \lambda^\star$.

\noindent \textbf{Lemma 6.} \textit{For any arbitrary ($\rho>1,\lambda^{(0)}>0$), there exist $T$ such that at iteration $t=T$, $\lambda^{(t)} > \lambda ^ \star$.}

\noindent \textit{Proof}. Depending on the initialization of $\lambda^{(0)}$, we have two cases:\\
\noindent 1) $\lambda^{(0)} > \lambda^{\star}$. In this case, at $T=1$ we have $\lambda^{(T)} > \lambda ^ \star$. 

\noindent 2) $\lambda^{(0)} \leq \lambda^{\star}$. In this case, if we select any $\rho > 1$, condition $\lambda^{(T)} > \lambda ^ \star$ is satisfied for $T = \ceil{ \frac{\log(\lambda^\star/\lambda^{(0)})}{\log \rho}}$, where $\ceil{.}$ is the ceiling operator.

\noindent Based on above arguments, the value of $T$ depends on the value of $\lambda ^ \star$ which is unknown. That is the reason we iteratively update $\lambda$ instead of choosing a fixed value.
\hfill\(\blacksquare\) \\

According to the above lemma, the choice of ($\rho>1,\lambda^{(0)}>0$) does not affect the final solution. Thus, there is no need to find their values via cross-validation. 
    
	\ifCLASSOPTIONcaptionsoff
	\newpage
	\fi
	\bibliographystyle{IEEEtran}
	\bibliography{Citations}
	
\end{document}